\documentclass[a4paper,11pt]{article}

\usepackage{amssymb,epsfig}
\usepackage{amsmath}
\usepackage{amscd}
\usepackage{amsthm}

\newtheorem{thm}{Theorem}[section]
\newtheorem{lem}[thm]{Lemma}
\newtheorem{prop}[thm]{Proposition}
\newtheorem{cor}[thm]{Corollary}

\theoremstyle{definition}
\newtheorem{defn}[thm]{Definition}
\newtheorem{remark}[thm]{Remark}

\newcommand{\pf }{{\bf Proof. }}

\newcommand{\ind}{\operatorname{ind}}
\newcommand{\coker}{\operatorname{coker}}
\newcommand{\Diff}{\operatorname{Diff}}
\newcommand{\End}{\operatorname{End}}
\newcommand{\Img}{\operatorname{Im}}

\newcommand{\supp}{\operatorname{supp}}

\newcommand{\RR}{\mathbb{R}}
\newcommand{\CC}{\mathbb{C}}

\newcommand{\E}{\mbox{$\cal E$}}
\newcommand{\F}{\mbox{$\cal F$}}
\newcommand{\G}{\mbox{$\cal G$}}
\newcommand{\M}{\mbox{$\cal M$}}
\newcommand{\N}{\mbox{$\cal N$}}
\newcommand{\X}{\mbox{$\cal X$}}
\newcommand{\W}{\mbox{$\cal W$}}
\newcommand{\Y}{\mbox{$\cal Y$}}

\title{
Cobordism invariance\\
of the family index\footnote{Research supported in part by Funda\c{c}\~ao para a Ci\^encia e Tecnologia through the grant FCT POCI/MAT/55958/2004. \textit{2000 Mathematics Subject Classification}: Primary 19K56, Secondary 58J20, 14D99. 
\textit{Keywords and phrases}: Family index. Cobordism.}
}
\author{Catarina Carvalho, \\
Department of Mathematics, 
Instituto Superior T\'ecnico, \\
Technical University of Lisbon
\\ email: {\tt ccarv@math.ist.utl.pt}}

\date{}

\begin{document}

\maketitle

\begin{abstract}
We give a $K$-theory proof of the invariance under cobordism of the family index. We consider elliptic pseudodifferential families on a continuous fibre bundle with smooth fibres $M\hookrightarrow \M \to B$, and define a notion of cobordant families using $K^1$-groups on fibrations with boundary. We show that the index of two such families is the same using properties of the push-forward map in $K$-theory to reduce it to families on $B\times \mathbb{R}^n$. 
 \end{abstract}

\section*{Introduction}\label{s int}

The invariance under cobordism of the Fredholm index has been a useful tool in index theory, both as a means to obtain index formulas
 and as an important step towards so-called relative index theorems.
 The particular case of twisted signature operators was crucial in Atiyah and Singer's first proof of the index theorem on closed manifolds \cite{a-s-cobd}. There are now many proofs of cobordism invariance for Dirac operators on closed manifolds, see for instance \cite{braverman-cobdsm,higson-cobdsm, moroianu-cobdsm, nicolaescu-cobdsm}, mostly relying on the geometric structure of the Dirac operator and not easily extended to other settings.
 In \cite{cc cbdsm}, we  gave a proof of cobordism invariance that applies to general elliptic pseudodifferential operators, under suitable conditions on their $K$-theory symbol classes defining what we called symbol-cobordism. See also  \cite{moroianu cc cbdsm} for an analytic formulation of this result using the calculus of cusp pseudodifferential operators on manifolds with boundary.

  In this paper, we give a $K$-theory proof of cobordism invariance for families of elliptic pseudodifferential operators on closed manifolds along the lines of \cite{cc cbdsm}.  In particular, we establish conditions on the symbol of a given elliptic pseudodifferential family on a boundary that yield the vanishing of its index in the $K$-theory of the base. The crucial point is to use push-forward maps and functoriality of the family index in $K$-theory, as proved by Atiyah and Singer as the main tool in their proof of the index theorem for families \cite{a-sIV}. 
We consider continuous families, in the spirit of \cite{a-sIV}, so we will be working in the setting of continuous fibrations with smooth fibre diffeomorphic to some given  manifold.
    We need  a well-behaved notion of boundary over the base space, in order to  define symbol-cobordism for families, and 
 it will be important to construct boundary preserving embeddings into Euclidean space over the base. One can then use the induced push-forward maps to reduce the problem to this setting. 
 Cobordism invariance for families is then a consequence of the functoriality of the index map and, moreover,  it depends only on properties of the symbol $K$-theory class of a given elliptic family on a boundary.
Moreover, we refine here the result in \cite{cc cbdsm}, in that we give also a condition at the operator level taking symbols of supended operators on manifolds with boundary
 
 Index theory for families has been a subject of recent study in the context of singular spaces, in particular, in the context of operators on fibered manifolds
 and, more generally, of pseudodifferential calculi on groupoids \cite{lmn-ctns fam grpds, melrose-piazza fam dirac, melrose-piazza fam dirac odd, melrose-rochon fam ind bdry, melrose-rochon fam cusp,  monthubert-grpds, nistor fam inv bdle lie grps, nwx, paterson-ctns fam grpds}.
In what regards cobordism invariance of family indices, there is an early result on families of Dirac operators used by Shih \cite{shih} to yield a weaker version of the index theorem for families (extending Atiyah and Singer's early cobordism proof).
More recently, cobordism invariance has been a key tool to study index theory on manifolds with boundary. Melrose and Piazza gave a proof of cobordism invariance for smooth families of Dirac operators on a boundary, both for the odd and even cases  \cite{melrose-piazza fam dirac, melrose-piazza fam dirac odd}, and used it to obtain an index theorem for families of Dirac operators on a manifold with boundary in the context of the $b$-calculus. 
Even if their proof uses Atiyah and Singer's index formula, in that they check that the topological index is cobordism invariant, the tools used are similar to the ones we apply here. They use functoriality of the topological index with respect to push-forwards (which reduces to functoriality of the Thom isomorphism) and, in the even case, the vanishing of the index relies on extension properties of the symbol class of the Dirac family on a boundary with respect to the $K_1$-symbol for the self-adjoint Dirac family on the fibered manifold with boundary. This extension condition on the symbols can be seen to be equivalent to the one established here for general elliptic pseudodifferential families (see \S\ref{s cbdsm inv}).  

An analytic formulation of cobordism invariance for pseudodifferential families can also be found in \cite{melrose-rochon fam ind bdry}, in the context of existence of invertible perturbations of cusp pseudodifferential families on manifolds with boundary. 
On a very different line, Hilsum \cite{hilsum-cbdsm fol} used the general framework of Hilbert modules over $C^*$-algebras to show cobordism invariance of the index of Dirac operators on foliated manifolds, that is, of continuous families of Dirac operators on the leaves. In this case, the so-called longitudinal index is an element of $K_0(C^*G)$, where $G$ is the holonomy groupoid and $C^*G$ is the groupoid $C^*$-algebra.

  The approach to cobordism invariance followed here highlights the fact that the vanishing of the index for an elliptic family on a boundary depends only on the existence of suitable extensions of the symbol class to the boundary. It is close in spirit with Atiyah and Singer's $K$-theory proof of the index formula, in that it relies on functoriality properties of the analytic index with respect to push-forwards and on the construction of suitable embeddings into Euclidean space. In that respect, it has the advantage of having immediate generalizations, namely to the non-compact and equivariant case. Moreover,  it would be of interest to determine if  one can relate it with the functoriality properties of the foliation index, which were crucial in the proof of the index theorem for foliations, in order to find $K$-theoretical conditions ensuring the invariance of the longitudinal index of pseudodifferential operators under foliation cobordism.\\
 
Given a closed smooth manifold $M$ and a compact base space $B$, we consider a manifold $\M$ over $B$ as a fibre bundle with fibre $M$ and structure group $\Diff(M)$; this is a family of smooth manifolds diffeomorphic to $M$ such that the smooth structure varies continuously over $B$. A pseudodifferential operator on $\M$ is a continuous family $P=(P_b)_{b\in B}$ of pseudodifferential operators on the fibres $M_b$ (Definition \ref{def fam psos}). If each $P_b$ is elliptic, then there is a well-defined index, $$\ind(P):= [\ker\tilde{P}]-[\CC^k]\in K^0(B),$$ where $\tilde{P}$ is a suitable perturbation of $P$ (see Definition \ref{defn ind}).
  Moreover, $\ind(P)$ only depends on the symbol class of $P$, $\sigma(P)\in K^0(T\M)$, where $T\M$ denotes the tangent space along the fibres.

 Let now $\X$ be a manifold over $B$, with fiber a manifold $X$ with boundary, $M=\partial X$ (we assume that the structure group of $\X$ preserves $M$). Then there is a correspondent continuous bundle of boundaries $\M$ with fiber $M$ and structure group $\Diff(M)$ (more precisely, the closed subgroup of those diffeomorphisms that extend to $X$). 
 We call $\M$ the \textit{boundary over $B$} of $\X$, $\M=\partial_B \X$. Of course, if the total space of $\X$ is smooth, then $\partial_B\X=\partial \X$. Many results on manifolds with boundary carry over to the boundary over $B$, namely that $\partial_B\X$ always has a collar in $\X$ and, in particular, 
 that $T\X_{\scriptsize{\M}}\cong T\M\times \RR$, where $T\X_{\scriptsize{\M}}$ is the restriction of the tangent space along the fibres of $\X$ to $\M$. 
 The constructions of \cite{cc cbdsm} can therefore be generalized to families and we obtain a map of restriction of symbols
 $$ u_{\scriptsize\X}: K^1(T\X)\to K^0(T\M)$$
 defined using maps of restriction to the boundary and the Bott isomorphism.
 Our main result states that:
 
 \begin{thm} \label{cdm inv int} Let $P$ be an elliptic family of pseudodifferential operators on a manifold $\M$ over $B$, with symbol $\sigma(P)\in K^0(T\M)$. If $\X$ is such that $\partial_B\X=\M$ and  $\sigma(P)\in \Img u_{\scriptsize\X}$, then $\ind(P)=0$.
 \end{thm}
  We say that the pair $(M,\sigma(P))$ as above is symbol-cobordant to zero. 
  Defining symbol-cobordant families in the obvious way, we have then that the family index is invariant under such relation. 
  
  To prove Theorem \ref{cdm inv int}, we show that  symbol-cobordism is invariant with respect to push-forward, as long as we consider boundary-preserving embeddings.
   We see that there always exist such embeddings into Euclidean space over $B$, and then that the relevant $K$-group in this case is zero. Functoriality of the family index with respect to push-forward then yields the result.
    
    As a consequence of Theorem \ref{cdm inv int}, we also obtain a condition for cobordism invariance at the level of symbols of families on $\X$ (Corollary \ref{cor susp ops}). The key point is to identify elements of $K^1(T\X)$ with symbol classes of suspended families on $\X$ and noting that in this case restriction to the boundary yields the indicial operator.
       
    As we have seen in \cite{cc cbdsm}, 
 cobordism invariance holds also on non-compact manifolds, considering operators that are multiplication outside a compact and taking the closure of a suitable $*$-algebra. In particular, we saw that the functoriality of the index with respect to push-forwards can be extended to this class. 
 (Note that in this case symbol cobordism is really a condition on the symbols, since any manifold is cobordant to zero through a non-compact cobordism.)
  Even though we do not pursue the non-compact approach in detail here, one can check that the compatibility of the index with push-forwards  given in \cite{cc cbdsm} can be extended to families of multiplication operators outside a compact, so that an analogue of  Theorem \ref{cdm inv int}  follows in this case as well
 (see Remark \ref{fam mult infty}). 
 Note that the closure of this class of operators can be used, through homotopy, to compute the index on large classes of operators and, moreover, it contains well-known pseudodifferential calculi.  
 See \cite{cc-thesis, cc-nistor} for details.

 It is important to mention that Moroianu obtained in \cite{moroianu cc cbdsm} 
 a result equivalent to the one in \cite{cc cbdsm} on closed manifolds, using quite different techniques. He used the calculus of cusp operators to show that an elliptic pseudodifferential operator on a boundary that has a suitable extension to a cusp pseudodifferential operator has zero index. This approach has the advantage of giving an explicit condition for cobordism invariance at the operator level. He also gave a $K$-theory formulation of this result and showed it is equivalent to the one given in \cite{cc cbdsm} in the closed case. There is a straightforward analogue for families of Moroianu's condition at the level of $K$-theory, and, using Theorem \ref{cdm inv int}, the same proof applies to show that it yields cobordism invariance for elliptic families.\\

\section{The index for families}\label{s an-ind}

We consider families of pseudodifferential operators as in \cite{a-sIV} (see also \cite{lawson-mich} for a detailed account). Let $B$ be a compact Hausdorff space, $M$ be a smooth compact manifold without boundary, $E$ a smooth vector bundle over $M$. We denote by $\Diff(M)$ the group of diffeomorphisms of $M$, endowed with the topology of uniform convergence on each derivative.
Also, $\Diff(M;E)$ denotes the subgroup of $\Diff(E)$ of those diffeomorphisms that carry fibres to fibres linearly;  choosing a connection on $E$, $\Diff(M;E)$ is a topological group.

\begin{defn}\label{def 1}
A \textit{manifold over} $B$ is is a fibre bundle $M\hookrightarrow\M\to B$ with fibre $M$ and structure group $\Diff(M)$.
A vector bundle $\E$ over $\M$ is said to be
  a \textit{smooth vector bundle (along the fibres)} if $\E$ defines a fibre bundle  $E\hookrightarrow\E \to B$ with fibre $E$ and structure group $\Diff(E;M)$.
\end{defn}
A manifold over $B$ is then a family of manifolds diffeomorphic to $M$ such that the smooth structure varies continuously. 
 If $\M=B\times M$, we call it a trivial family; locally, every manifold over $B$ is of this form. 
Note also that a smooth vector bundle is a continuous family of smooth vector bundles over $M$; moreover, given such a bundle $E\hookrightarrow\E \to B$, the map $\Diff(E;M)\to \Diff(M)$ induces a manifold $\M$ over $B$.  
The cotangent and tangent bundles along the fibres, denoted by $T^*\M$, $T\M$, respectively, are smooth vector bundles over $\M$, as in Definition \ref{def 1}.

A {submanifold} of a manifold $M\hookrightarrow\M\to B$ is just a sub-bundle. It is clear that if $N\subset M$ is closed, then there is a submanifold $\N$ of $\M$ with fibre $N$ if, and only if, the structure group of $\M$ can be reduced to the closed subgroup $\Diff(M,N)$ of diffeomorphisms that preserve $N$. 
 
Let  now $E$, $F$ be vector bundles over $M$, with $ \Gamma(E)$, $\Gamma(F)$ the spaces of smooth sections of $E$, $F$. We denote by  $\Psi^m(M;E,F)$ the space of order $m$ classical pseudodifferential operators $P: \Gamma(E)\to \Gamma(F)$, as in \cite{horm-book} .
It is a Fr\'echet space with the topology induced by the semi-norms of local symbols and their derivatives
$$\|P\|_{K, \alpha, \beta}:= \left |\frac{\partial_x^\alpha \partial_\xi^\beta p(x,\xi)}{(1+|\xi|)^{m-|\beta|}}\right|, $$
where $p$ is the classical symbol inducing $P$ on some coordinate chart $U$ of $M$ trivializing $E$ and $F$,  $K\subset U$ is compact, and for any multi-indices $\alpha$, $\beta$.

If we let $\Diff(E,F;M)$ be the subgroup of $\Diff(E\oplus F;M)$ of those diffeomorphisms that map $E$ to $E$ and $F$ to $F$, then $\Diff(E,F;M)$ acts on $\Psi^m(M;E,F)$ by $P\mapsto f_1^{-1}Pf_2$, for $f=(f_1, f_2)$  and the action is continuous (see \cite{a-sIV}). 
 To each pair of continuous families $\E$, $\F$ of vector bundles over $M$, with fibres diffeomorphic to $E$, $F$, we can then associate a fibre bundle $\Psi^m(\M;\E,\F)\to B$ with fibre $\Psi^m(M;E,F)$ and structure group $\Diff(E,F;M)$. The manifold $\M$ over $B$ is induced by the map $\Diff(E,F;M)\to \Diff(M)$.
 
 \begin{defn}\label{def fam psos}
A \textit{continuous family of pseudodifferential operators} on $\M$ is a continuous section $P$ of $\Psi^m(\M;\E,\F)$; we write $P=(P_b)_{b\in B}$. The family $P$ is said to be \textit{elliptic} if each $P_b$, $b\in B$, is elliptic.
\end{defn}

If $\M=B\times M$, $\E=B\times E$, $\F=B\times F$, then a continuous family is just a continuous map $P: B\to{\Psi}^m(M;E,F)$. Every family is locally of this form.
 
 Recall also from \cite{horm-book} that there is a surjective symbol map 
 $$\sigma: \Psi^m(M;E,F) \to S^m(M;E,F),$$ 
 where $S^m(M;E,F)$ is the space of classical symbols, that is, smooth maps on the cotangent bundle $T^*M\backslash 0$ with values in $\End(E,F)$, that are positively homogeneous of degree $m$. The map $\sigma$ is invariant under the action of $\Diff(E,F;M)$
  and the action is continuous (endowing $S^m(M;E,F)$ with the $\sup$-norm topology on the sphere bundle $SM$),
 so that we get a fibre bundle $S^m(\M;\E,\F)$ over $B$ with fibre $S^m(M;E,F)$ and group $\Diff(E,F;M)$.

 \begin{defn}\label{defn symbol}
 The \textit{symbol} $\sigma_B(P)$ of a family $P=(P_b)_{b\in B}$ is the continuous section of $S^m(\M;\E,\F)$ given fibrewise by $\sigma(P_b)$.
\end{defn}

The family of smooth symbols  $\sigma(P_b)$  depends continuously on $b\in B$.  Elements of  $S^m(\M;\E,\F)$ are maps $T^*\M\to \End(\pi^*\E,\pi^*\F)$, with $T^*\M$ the cotangent bundle along the fibres and $\pi: T^*\M\to \M$ the projection. If the family is elliptic, then each $\sigma(P_b)$ is invertible outside the zero-section, so that $\sigma_B(P)$ is invertible outside a compact in $T^*\M$.
Moreover, endowing $\M$ with a metric, that is, with a continuous family of metrics on the fibres $TM$, we can identify $T^*\M$ with $T\M$ and (elliptic) symbols reduce to (invertible) maps on the sphere bundle $S\M$.

Finally, we now define the index of an elliptic family. Let $P$ be an elliptic family in $\Psi^m(\M; \E,\F)$,
 so that each $P_b$, $b\in B$, is Fredholm. If $\dim \ker P_b$ were locally constant, then the family of vector spaces $\ker P_b$ would form a vector bundle over $B$, and the same for $\coker P_b$; in this case, one could define  the index of the family $P$ as the $K$-theory class $[\ker P_b]-[\coker P_b]\in K^0(B)$. 
 In the general case, one can define an elliptic family $\tilde{P}: \Gamma(\E)\oplus \CC^k\to \Gamma(\F)$ such that  $\tilde{P}_b$ is surjective, for each $b\in B$, where
\begin{equation}
\label{tildeP}
\tilde{P}_b(u; \lambda_1,...,\lambda_k):= P_b(u)+\lambda_1w_1(b)+...+ \lambda_kw_k(b),
\end{equation}
for some smooth sections $w_i$ of $\F$, $i=1,...,k$.
  Since each $\tilde{P}_b$ is surjective,  the family $\ker \tilde{P}_b$ indeed defines a vector bundle over $B$ (each $\tilde{P}_b$ is Fredholm and the Fredholm index is locally constant). Moreover, the class $[\ker \tilde{P}]-[B\times \CC^k ]\in K^0(B)$ depends only on $P$ (see \cite{a-sIV, lawson-mich}). 

\begin{defn}\label{defn ind}
The \textit{analytic index of the elliptic family} $P$ is given by
\begin{equation*}
\ind (P):=[\ker \tilde{P}]-[B\times \CC^k ] \in K^0(B),
\end{equation*}
where $\tilde{P}: \Gamma(\E)\oplus \CC^k\to \Gamma(\F)$ is as in (\ref{tildeP}).
\end{defn}
The index defined above is homotopy invariant, hence it depends only on the (homotopy class) of the symbol $\sigma_B(P)$. 
 Moreover, it can be seen that  the index only depends on the values of the symbol on $S\M$, that is,  the index of a family $P$ does not depend on the order of $P$. Therefore, as far as computing indices goes, it is just as good taking operators of order $0$.

Now, the symbol of an elliptic family  $\sigma_B(P)$ is invertible outside a compact, hence it defines a class 
\begin{equation}\label {symb class}
[\sigma_B(P)]:=[\pi^*\E, \pi^*\F, \sigma_B(P)]\in K^0(T\M)
\end{equation}
(where we identify $T^*\M$ with $T\M$). Since $ K^0(T\M)$ is exhausted by symbol classes, just as in the operator case, 
we have then a well-defined family index map
\begin{equation}\label {ind map family}
\ind : K^0(T\M)\to K^0(B),\quad [\pi^*\E, \pi^*\F, \sigma_B(P)]\mapsto \ind(P),
\end{equation}
which is, moreover, a group homomorphism. 
We will show in this paper that  the family index is cobordism invariant under a suitable notion of cobordant families that depends only on the manifolds involved and on the  $K$-theory class of the family symbol. To do this, we will make use of push-forward maps.

\section{Push-forward maps}\label{s i!}
We consider now a manifold $\X$ over $B$ with fiber some smooth manifold $X$ \textit{possibly with boundary}. In this case, we assume that the structure group of $\X$ is $\Diff(X,\partial X)$, the subgroup of $\Diff(X)$ of those diffeomorphisms that leave $\partial X$ invariant.
 Also, the structure group of a smooth bundle $\E$ over $\X$ is given in this case by $\Diff(E,X; \partial E)$ of those elements of $\Diff(E,X)$ that preserve $\partial E=E_{|\partial X}$.

Recall that an embedding $i: X\to Y$ is neat, following \cite{hirsch}, if $\partial X= X\cap \partial Y$ and if $X$ intersects $\partial Y$ transversally; in this case, the boundaries are also embedded and $X$ always has an open tubular neighborhood in $Y$. 

Let  $i: \X\to B\times V$ be a morphism,
where $V$ denotes either $\RR^m$ or $\RR^{m}_+$, the positive half-space in $\RR^m$, in case $X$ has a boundary. We say that $i$ is a \textit{neat embedding of manifolds over $B$} if $i$ restricts to neat embeddings $i_b: X_b\to V$ on each fiber of $\X$. Note that in this case, $i$ is a homeomorphism onto $i(\X)$, smooth on each fiber.

  We identify here the open tubular neighbourhood of $X_b$ in $V$ with the normal bundle $N_b=N$. We obtain then a continuous family $\N$ of vector bundles over $X$ with fibre $N$; moreover, $\N$ is open in $B\times V$.

The induced embedding $T\X \to B\times TV$ is such that the fiberwise embeddings are also neat. The normal bundle of $TX_b$ in $TV$ is given by $TN_b\cong \pi_X^*(N_b\oplus N_b)$, where $\pi_X: TX\to X$ is the projection, so that the vector bundle $T\N$ can be identified with $\pi_{\scriptsize\X}^*(\N\oplus \N)$, $\pi_{\scriptsize\X}: T\X\to \X$, and therefore it has a complex structure.
We have then  a Thom isomorphism 
\begin{equation}
\label{thom iso}
\rho: K^0(T\X)\to K^0(T\N).
\end{equation}
Composing with the $K$-theory map  induced by the open inclusion $T\N\to B\times TV$
\begin{equation}
\label{ext map}
h: K^0(T\N)\to K^0(B\times TV),
\end{equation}
one obtains the \textit{push-forward map}
\begin{equation}\label{i!}
i_!=h\circ \rho : K^0(T\X)\to K^0(B\times TV).
\end{equation}
One can also easily define  a push-forward map between $K^1$-groups: simply take the induced embedding $T\X\times \RR \to B\times TV\times\RR$ to get $i_!: K^1(T\X)\to K^1(B\times TV)$.

All the constructions above work in the same way for maps $i: \X \to \Y$, for some manifold $\Y$ over $B$, as long as $i$ restricts to neat embeddings on the fibers. Here we only need the case $\Y= B\times V$; we will see later (Corollary \ref{emb}) that such an embedding into a trivial family of manifolds always exists,
 and in fact it can be chosen so as to yield an embedding $\M\to \RR^n$ of the boundaries over $B$.

One of the crucial properties of the push-forward map is its functoriality with respect to the index.
In their proof of the index theorem for families \cite{a-sIV}, Atiyah and Singer show that:
 
\begin{thm}[Atiyah-Singer]\label{i! vs ind}
Let $\M$ be a manifold over $B$, with fibre $M$ compact without boundary. Then the following diagram commutes
\begin{equation}
\begin{CD}
K^0(T\M)         @>{i_!}>>   K^0(B\times T\RR^m)          \\
@V{\ind}VV            @VV{\ind}V           \\
K^0(B)           @=           K^0(B)     
\end{CD}
\end{equation}
where $\ind$ stands for the analytic index for families, as in (\ref{ind map family}).  
\end{thm}
The  map on the right-hand side computes the index of families of operators that are the identity outside a compact of $\RR^m$ (see also Remark \ref{fam mult infty} below). The families index theorem in this setting states that it coincides with the Bott isomorphism.
The proof of Theorem \ref{i! vs ind} goes much as in the case of the classic index theorem, taking extra care with the added sections in the definition of the family index.

\begin{remark}[Operators that are multiplication at infinity]\label{fam mult infty}

	Let $M$ be a non-compact manifold, $E$, $F$ vector bundles on $M$.
 In \cite{cc cbdsm}, we considered the class $\Psi_{mult}(M;E,F)$ of pseudodifferential operators on $M$ that are multiplication at infinity. For an operator $P: C_c^\infty(M;E) \to C^\infty(M;F)$, we have $P\in \Psi_{mult}(M;E,F)$ if, and only if, there is a compact $K\subset M$ and $p: E\to F$ such that 
\begin{equation}
\label{psos mult infty}
P=P_1 + p,
\end{equation}
where $P_1: C_c^\infty(M;E) \to C^\infty(M;F)$ is a pseudodifferential operator on $M$ such $\Img P_1\subset C_K^\infty(M;F)$, $\Img P_1^*\subset C_K^\infty(M;E)$, that is, $P_1$ has compactly supported Schwartz kernel. This means that if $\supp u\subset M\setminus K$, then $Pu=pu$, $P^*u=p^*u$.
In particular, $\sigma(P)$ is independent of $\xi$ for $x\notin K$.

Let now $\M$, $\E$, $\F$ be as before, with non-compact $M$, over some compact base space $B$, and consider the class $\Psi_{mult}(\M;\E,\F)$ of continuous families of pseudodifferential operators that are multiplication at infinity, with bounded symbols. Since $B$ is compact, given such a family $P=(P_b)$, we can pick a compact $K\subset M_b$, independent of $b\in B$, in the above definition.
We say that $P\in \Psi_{mult}(\M;\E,\F)$ is \textit{fully elliptic} 
 if, for each $ b\in B$, $\sigma(P_b)$ is invertible on $T^*M\setminus 0$ and the inverse is bounded, which means that it also independent of $\xi$ for $x\notin K$. Since $\sigma(P_b)$ is constant on fibers over $M\setminus K$, it is invertible outside the (fixed) compact $\pi^{-1}(K)\subset T^*M$, with $\pi : B^*M\to M$ the ball cotangent bundle.

Hence, if $P$ is a fully elliptic family, then the symbol $\sigma_B(P)$ is invertible outside a compact in $T\M$ and it defines a $K$-theory class. Moreover, the family index (as in Definition \ref{defn ind}) is well defined, since we have in this case $\ker P_b \subset  C_K^\infty(M;E)$, $\coker P_b  \subset  C_K^\infty(M;F)$, and it depends only on $\sigma(P)$. 

The excision property of the index given in \cite{a-sIV} holds in this setting and one can use it to extend  the compatibility of the family index with push-forwards, exactly as in \cite{cc cbdsm} (Theorem 3.9), since all the constructions carry over nicely to families. We obtain then, using the notation of Theorem \ref{i! vs ind}, with $M$ not necessarily compact,
\begin{equation}
\label{i! vs ind mult infty}
\ind \circ i_! = \ind
\end{equation}
where now we are computing the index of families of operators that are multiplication at infinity on $\M$.
\end{remark}

\section{Cobordant families}\label{s cbdt fam}

We define here the notion of symbol-cobordism for families of pseudodifferential operators, generalizing \cite{cc cbdsm}; one of its main features is that it is preserved by push-forward maps.

Let $X$ be a manifold with boundary and $M=\partial X$. If $\X$ is a manifold over $B$ with fibre $X$, the inclusion $M\to X$ induces a submanifold $\M$ over $B$, with fibre $M$ and structure group $\Diff(M)$. 
We call $\M$  the \textit{boundary over $B$} of $\X$, and write $\M=\partial_B \X$. 
If $X\hookrightarrow \X\to B$ is given by a smooth fibration, then $\partial_B\X= \partial \X$, the induced boundary fibration. 

Note that if $\M=\partial_B \X$, then the structure group of $\M$ can be reduced to the closed subgroup of diffeomorphisms that can be extended to $X$. 
Conversely, given $M$ and $X$ such that $\partial X=M$, then there exists $\X$ with $\partial_B \X=\M$ if, and only if, the structure group of $\M$ can be reduced to the subgroup of  diffeomophisms that can be extended to $X$.

Before proceeding to give our main definition, we show that a continuous family of boundaries always has a collar neighborhood over $B$.

\begin{prop}\label{collar nbhd over B}
Let $\M$, $\X$ be manifolds over $B$, with fibres $M$, $X$, such that $\M=\partial_B\X$. Then, there is an isomorphism 
$$\Phi: \M\times [0,1)\to \W\subset \X,$$
where $\W$ is a submanifold of $\X$ with fibre $W\subset X$, open, and structure group $\Diff(\overline{W},M)$, such that $\Phi$ restricts to a collar neighborhood of the boundary on each fibre.
\end{prop}

\pf 
Let ${\phi}: M\times [0,1] \to \overline{W}\subset X$ be a collar neighborhood of $M$ in $X$. We first note that the structure group $\Diff(X,M)$ of $\X$ can be reduced to $\Diff(X,\overline{W})$: for each $\varphi\in \Diff(X,M)$, it is easy to check that there are $0<\alpha, \beta\leq 1$ such that $\varphi(W_\alpha) = W_\beta$, with $W_\alpha \cong M\times [0,\alpha]$, $W_\beta \cong M\times [0,\beta]$; pick $\psi_\alpha, \psi_\beta \in  \Diff(X,M)$ such that $\psi_\alpha(\overline{W})=W_\alpha$, $\psi_\beta(\overline{W})=W_\beta$. We have then $\psi_\beta^{-1}\varphi\psi_\alpha\in \Diff(X,\overline{W})$. 

We conclude that there is a submanifold $\overline{\W}$ of $\X$, with fibre $\overline{W}$.
To show that there is a globally defined map ${\Phi}: \M\times [0,1]\to \overline{\W}\subset \X$, we need to check compatibility with transition functions. Note that the transition functions for $\M$ and $\overline{\W}$ are given by the restriction of transition functions of $\X$, so that if $\varphi\in \Diff(X,\overline{W})$ is a transition function for $\X$, we want to see that
\begin{equation}
\label{ }
\varphi\circ {\phi}=  {\phi}\circ (\varphi_{|M},1).
\end{equation}
Noting that $\varphi^{-1}\circ {\phi}\circ (\varphi_{|M},1): M\times [0,1]\to \overline{W}$ is also a collar for $M$ in $X$, it follows from \cite{hirsch} (Theorem 8.1.8) 
 that there is $\varphi^\prime$ diffeotopic to $\varphi$, with $\varphi^\prime=\varphi$ on $M$ and on $X\setminus W$, such that 
\begin{equation}
\label{qqq}
{\phi}=( \varphi^\prime)^{-1} {\phi}\circ (\varphi^\prime_{|M},1)
\end{equation}
on $M\times [0,\alpha]$, for some $0<\alpha\leq 1$. Identifying $M\times [0,\alpha]$ with $M\times [0,1]$ through diffeomorphisms of $X$, we can assume that (\ref{qqq}) holds on $\M\times [0,1]$.
 
 Since the local smooth structure of $\X$ remains invariant under diffeotopies, it follows that there is a map
$ {\Phi}: \M\times [0,1] \to \overline{\W},$
with $\overline{\W}$ a manifold over $B$ with fibre $\overline{W}$ and structure group $\Diff(\overline{W},M)$, such that ${\Phi}$ restricts to ${\phi}$ on each fibre. Moreover, $\Phi$ is an isomorphism. 
 Restricting ${\Phi}$ to $\M\times [0,1)$ we get the result.
\qed

Let then $\M$ and $\X$ be manifolds over $B$ with $\M=\partial_B \X$. It is easy to check that $T\X_{\scriptsize{\M}}$, the restriction of the vertical tangent bundle $T\X$ to $\M$, is also a sub-bundle of $T\X$, with fiber $TX_{|M}$ (and structure group $\Diff(TX_{|M},M)$).
The inclusion $T\X_{\scriptsize\M}\to T\X$ induces a map of restriction to the boundary in $K$-theory
\begin{equation}\label{rest map}
r_{\scriptsize{\X,\M}}: K^j(T\X)\to K^j(T\X_{\scriptsize{\M}}).
\end{equation}
As in \cite{cc cbdsm}, it is not hard to show that push-forward maps behave well with respect to restriction to the boundary.

\begin{prop}\label{i! vs rest}
Let $\M=\partial_B \X$. If $i: \X \to B\times V$ is an embedding of manifolds over $B$ that restricts to an embedding $\M\to B\times W$, where $W=\partial  V$, then the following diagram commutes
\begin{equation*}
\begin{CD}
K^j(T\X)                   @>{i_!}>>   K^j(B\times TV)          \\
@V{r_{\scriptsize{\X,\M}}}VV                 @VV{r_{V,W}}V           \\
K^j(T\X_{\scriptsize{\M}})        @>{i_!}>>           K^j(B\times TV_W).    
\end{CD}
\end{equation*}
\end{prop}

\pf  Functoriality of the Thom isomorphism yields that it commutes with restriction maps. To see that the same happens with respect to the $K$-theory maps $h_X: K^j(T\N)\to K^j(B\times TV)$ and $h_M:K^j(T\N_{\M})\to K^j(B\times TV_W)$, as in  (\ref{ext map}), we just check that $r_{V,W}\circ h_X$ and $h_M\circ r_{\scriptsize{\X,\M}}$ can be written as maps in $K$-theory induced by the same morphism.
\qed 

Now, it is well known that there is a smooth isomorphism $TX_{|M}\cong TM \times \RR$, for $M=\partial X$. It follows from Proposition \ref{collar nbhd over B}
 that the same holds for tangent bundles over $B$: if $\Phi: \M\times [0,1) \to \W$ is a collar neighborhood over $B$, $\W\subset \X$ with open fiber $W$, then
 $$ T\X_{\scriptsize{\M}}=T\W_{\scriptsize{\M}}\cong T(\M\times [0,1))_{\scriptsize{\M}\times \{0\}}=T\M \times \RR,$$
 where $T\M\times \RR$ is the vector bundle over $\M$ given by $T\M\oplus (M\times \RR)$, so that the fiber is $TM\times \RR$ and the structure group 
 acts trivially on $\RR$. 
%
The restriction map  becomes now
$$ r_{\scriptsize{\X},\scriptsize{\M}}: K^1(T\X)\to K^1(T\M\times \RR)=K^0(T\M\times\RR^2).$$
Taking the Bott isomorphism $\beta: K^0(T\M)\to K^0(T\M\times\RR^2)$, we consider the map
\begin{equation}\label{map rest symb}
u_{\scriptsize{\X, \M}}:=\beta^{-1} r_{\scriptsize{\X,\M}}: K^1(T\X)\to K^0(T\M),
\end{equation}
which in \cite{cc cbdsm} was referred to as restriction of symbols. We will see later that elements of $K^1(T\X)$ can be regarded as symbol classes of suspended operators. 
 
 Let now, for $i=1,2$, $M_i$ be a manifold without boundary and $\M_i$ be a manifold over $B$ with fibers $M_i$. If $\X$ is such that $\partial_B \X=\M_1 \sqcup \M_2$, then the structure group of $\X$ can be reduced to the subgroup of diffeomorphisms that leave both $M_1$ and $M_2$ invariant, so that $\M_1$, $\M_2$ are submanifolds of $\X$, and it is easy to check that there are well-defined maps of restriction of symbols
 \begin{equation}\label{map rest symb 1,2}
u_{\scriptsize{\X, \M_i}}: K^1(T\X)\to K^0(T\M_i),
\end{equation}
for $i=1, 2$, and that $u_{\scriptsize{\X, \M}}=_{\scriptsize{\X, \M_1}}\oplus \; u_{\scriptsize{\X, \M_2}}$.
Let $P_i$  be an elliptic family of pseudodifferential operators on $\M_i$ with symbol $\sigma_i=\sigma(P_i)\in K^0(T\M_i)$, $i=0,1$.
 We give the following definition: 

\begin{defn}\label{def cbdt families} 
We say that $(\M_1, P_1)$ and $(\M_2, P_2)$ are \textit{symbol-cobordant} families if there is a manifold $\X$ over $B$, with fiber $X$, and $\omega \in K^1(T\X)$, such that

(i) $\partial_B \X = \M_1  \sqcup \M_2$ (in particular, $\partial X=M_1  \sqcup M_2$);

(ii) $\sigma_1=u_{\scriptsize{\X, \M_1}}(\omega)$, $\sigma_2= - \,u_{\scriptsize{\X, \M_2}}(\omega)$.
\end{defn}
We write $(\M_1,\sigma_1)\sim(\M_2, \sigma_2)$. Note that $(\M_1,P_1)$ and $(\M_2,P_2)$ are cobordant  if, and only if, $(P_1\oplus P_2^*, M_1  \sqcup M_2)$ is cobordant to zero, where $P_2^*$ is the adjoint family and $P_1\oplus P_2$ is a family on $\M_1  \sqcup \M_2$ (we have $-\sigma_2=\sigma(P_2^*)$). 
As in \cite{cc cbdsm}, the definition above defines an equivalence relation on pairs $(\M,P)$ within manifolds with the same fiber dimension and the collection of equivalence classes is an abelian group.

We now check that the notion of cobordant families is preserved by push-forward maps. 

\begin{prop}\label{i! vs u} Let $\X, \M$ be manifolds over $B$ with $\partial_B\X=\M$ and let $i: \X \to B\times \RR^{m+1}_+$ be an embedding of manifolds over $B$ that induces an embedding $i: \M\to B\times \RR^m$. Then the following diagram commutes
\[
\begin{CD}
K^1(T\X)                           @>{i_!}>>        K^1(B\times T\RR^{m+1}_+)     \\
@V{u_{\scriptsize{\X,\M}}}VV                                     @VV{u_m}V           \\
K^0(T\M)         @>>{i_!}>           K^0(B\times T\RR^m),   
\end{CD}
\]
where $u_{\scriptsize{\X,\M}}$ and $u_m$ are maps of restriction of symbols, as in (\ref{map rest symb}).
\end{prop}

\pf  It follows directly from (\ref{i! vs rest}) and the fact that push-forward maps commute with the Bott isomorphism.
\qed

\begin{cor}\label{cobdsm vs i!}
In the conditions of Proposition \ref{i! vs u}, if $\sigma\in K^0(T\M)$ and $(\M,\sigma)\sim 0$ then $(B\times\RR^m, i_!(\sigma))\sim 0$.
\end{cor}

\pf Straightforward from the definition of symbol-cobordism and Proposition \ref{i! vs u}.
\qed

\section{Cobordism invariance}\label{s cbdsm inv}

We now show that if a given elliptic family is symbol-cobordant to zero, then its index is zero. The main idea  is to use Corollary \ref{cobdsm vs i!} to reduce the proof to Euclidean space where it will be trivial. 

We start with showing that an embedding such as in Proposition \ref{i! vs u} indeed exists.
Let then $\X$ be a manifold over $B$ with $\partial_B \X=\M$, $\X$, $\M$ with fibres $X, M$, respectively. Also, let $i: X\to \RR^k$ be an embedding, which always exists by Whitney's theorem. It is easy to show  that, for compact $B$,  one can define an embedding
\begin{equation}
\label{ext emb}
\tilde{i}: \X\to B\times \RR^m,\qquad m\geq k,
\end{equation}
which restricts to an embedding on each fiber $X_b$ of $\X$ (see \cite{a-sIV}).
We need however an embedding that restricts to an embedding of the boundary $\M$ over $B$. In \cite{cc cbdsm}, we constructed such an embedding using a  collar neighborhood of $M$ in $X$, which is what we shall do here. We use the collar neighborhood over $B$ given in Proposition \ref{collar nbhd over B} to define  a map
\begin{equation}
\label{tildealpha}
\tilde{\alpha}: \X \to[0,1],\qquad \M=\alpha^{-1}(0),
\end{equation}
such that on each fiber, $\tilde{\alpha}$ restricts to a defining function of the boundary $\tilde\alpha_b=\alpha : X_b\to [0,1]$, with $\alpha^{-1}(0)=M$ and $d\alpha\neq 0$ on $M$. 

\begin{prop}\label{emb}
Let $\X$ be a manifold over $B$ with $\partial_B \X=\M$. There exists a neat embedding over $B$
\begin{equation*}
h: \X \to B\times \RR^{m+1}_+, 
\end{equation*}
which induces an embedding  $\M\to B\times \RR^m$. 
\end{prop}

\pf Let $\tilde{i}: \X\to B\times \RR^m$ be an embedding over $B$, as in (\ref{ext emb}) and $\tilde{\alpha}: \X \to[0,1)$ as in (\ref{tildealpha}). Define
$$h:=(\tilde{i}, \tilde{\alpha}): \X \to B\times \RR^{m+1}_+.$$
Then $h$ restricts to neat embeddings $h_b=(i,\alpha)$, where $i: X\to \RR^m$ is an embedding and $\alpha$ a boundary defining function. 
\qed

We use the embedding given above, together with Proposition \ref{i! vs u} to reduce the proof of cobordism invariance to $B\times \RR^m$. It is trivial in this case, as the following lemma will show.

\begin{lem} $K^0([0,1)\times U)=0$, for any locally compact space $U$.
\end{lem}
\pf 
We have $K^0([0,1)\times U)=K_0(C_0([0,1)\times U))$, where $C_0([0,1)\times U)$ is the $C^*$-algebra of continuous functions on $[0,1)\times U)$ that vanish outside a compact set. Since $C_0([0,1)\times U)\cong C_0([0,1),C_0(U))$, which is just the cone of the $C^*$-algebra $C_0(U)$ and hence, it is homotopy equivalent to $0$, 
  the result follows for $K^0$, and also for $K^1$, writing $K^1([0,1)\times U)=K^0([0,1)\times U\times\RR)=0$.
\qed

In particular, it follows from the lemma above that, for any (locally) compact base space $B$ and $m\in\mathbb{N}$, we have
\begin{equation}\label{}
K^1(B\times T\RR^{m+1}_+)=K^0(B\times \RR^{2m+1}\times [0,1))=0.
\end{equation}
We are now ready to prove:

\begin{thm}\label{thm cbdsm inv}
Let $\M$ be a manifold over $B$ and $P\in \Psi(\M;\E,\F)$ be an elliptic family of pseudodifferential operators over $B$ with symbol $\sigma=\sigma_B(P)\in K^0(T\M)$. If $(\M,\sigma)$ is symbol-cobordant to zero, then $\ind(P)=0$. 
\end{thm}

\pf From the definition of symbol-cobordism, we know that there is a manifold $\X$ over $B$, with fiber $X$, such that $\partial_B \X=\M$ 
and a class $\omega \in K^1(T\X)$ such that $\sigma=u_{\scriptsize{\X,\M}}(\omega)$. 

Considering an embedding $h: \X \to B\times \RR^{m+1}_+$, as in Corollary \ref{emb},
we have from Corollary \ref{cobdsm vs i!} that $(B\times \RR^m, h_!(\sigma))$ is also cobordant to zero, so that $h_!(\sigma)=u_m(h_!(\omega))$ (we are using the notations of Proposition \ref{i! vs u}). But $h_!(\omega)\in K^1(B\times T\RR^{m+1}_+)=0$, therefore $h_!(\sigma)=0$. Since the index is invariant under push-forward maps (Theorem \ref{i! vs ind}) we have finally
$$\ind(P)=\ind(\sigma)=\ind(h_!(\sigma))=0.$$
\qed

\begin{cor}
Let $P_1\in \Psi(\M_1;\E_1,\F_1)$, $P_2\in \Psi(\M_2;\E_2,\F_2)$ be elliptic families of pseudodifferential operators over $B$ with symbols $\sigma_1, \sigma_2$, respectively. If $(\M_1, \sigma_1)$ is symbol-cobordant to $(\M_2,\sigma_2)$, then $\ind(P_1)=\ind(P_2)$. 
\end{cor}

\pf It follows from $(P_1\oplus P_2^*, \M_1\sqcup \M_2)$ being cobordant to zero, so that $\ind(P_1\oplus P_2^*)=0$, and the additivity of the family index.
\qed

In the spirit of Remark \ref{fam mult infty}, we note that the definition of symbol-cobordism and the results given thereafter are purely $K$-theoretical and do not depend in any way on the compacity of the fibres $M$ of $\M$. Hence, Theorem \ref{thm cbdsm inv} holds for non-compact $\M$ in the setting of pseudodifferential operators that are multiplication at infinity and we state it here as a corollary. 

\begin{cor}
Let $\M$ be a manifold over $B$ with non-compact fibres and $P=(P_b)$ be a family of fully elliptic operators in $\Psi_{mult}(\M;\E,\F)$,with symbol $\sigma=\sigma_B(P)\in K^0(T\M)$. If $(\M, \sigma)\sim 0$, we have $\ind (P)=0$.
\end{cor}
In particular, we obtain straightaway  the invariance of the index for families that are homotopic to multiplication families outside a compact, for instance for families of scattering operators; see \cite{cc-nistor} for details.

Let  $X\hookrightarrow \X \to B$ be a manifold over $B$ ,with $X$ a manifold with boundary. We now identify classes in $K^1(T\X)$ with symbols of suspended operators. We work here in the setting of the $b$-calculus (see \cite{melrose APS} for an extended treatment). We assume that $\X$ is endowed with a continuous family of exact $b$-metrics and consider the class $\Psi_b(\X,\G)$ of families of $b$-pseudodifferential operators $\Psi_b(X,G)$ on the fibers, as in the boundaryless case. The symbol of such a family is defined on the $b$-tangent space along the fibers $^bT\X$; since there is a (non canonical) bundle isomorphism $^bT\X\cong T\X$, at the $K$-theory level we can work with $T\X$. 
To each operator $Q_b\in\Psi_b(X,G)$, $b\in B$, one can associate a translation invariant pseudodifferential operator on $\partial X\times \RR$, that is, a 1-parameter family of pseudodifferential operators on $M=\partial X$, the indicial family 
$I_M(Q_b)\in \Psi_{sus}(M, G_{M})$ (where the convolution kernels are assumed to vanish rapidly with all derivatives at infinity). It follows from the definition that there is a compatibility condition with respect to the principal symbol
\begin{equation}
\label{symb I(Q) vs Q}
\sigma(I_M(Q))=\sigma_B(Q)_{|\scriptsize{\M}}.
\end{equation}
In particular, if $Q$ is elliptic then $I_M(Q)$ also is and $[\sigma(I_M(Q))]\in K^0(T(\partial \X)\times \RR)$.
If we now take suspended elliptic families of $b$-operators on $\X$, that is, a family $Q=(Q_t)$ on $X\hookrightarrow \X^\prime=\X\times [0,1] \to B^\prime = B\times [0,1]$, depending smoothly on $t$, and assume that $Q_{0}, Q_1$ are invertible, then the symbol $\sigma_{B^\prime}(Q)$ defines a class in $K^0(T\X \times (0,1))=K^1(T\X)$. Moreover, $K^1(T\X)$ is exhausted by such symbol classes. Note that $I_M(Q)$ is now defined on the double suspension of $\M$.

\begin{cor}\label{cor susp ops}
Let $\M$ be a manifold over $B$ and $P\in \Psi(\M;\E,\F)$ be an elliptic family with symbol $[\sigma_B(P)]\in K^0(T\M)$. If $\M=\partial_B \X$, for some manifold $\X$ over $B$ and there is an elliptic family $Q\in \Psi_b(\X\times [0,1];\G)$, invertible at $t=0,1$, with indicial family $I_M(Q)$, such that
$$[\sigma(I_M(Q))]=[\sigma_B(P)] \cup \beta$$
with $\beta$ the Bott class, then $\ind(P)=0$. 
\end{cor}

\pf 
It follows from (\ref{symb I(Q) vs Q}) that $[\sigma(I_M(Q))]=r_{\scriptsize{\M}}[\sigma_{B^\prime}(Q)]$, with $r_{\scriptsize{\M}}$ the map of restriction to the boundary. Hence, $[\sigma_B(P)]= u_{\scriptsize{\M}}[\sigma_{B^\prime}(Q)]$, with $u_{\scriptsize{\M}}$ the symbol restriction map (\ref{map rest symb}), which gives that  $(\M, [\sigma_B(P)])$ is symbol-cobordant to zero.
\qed\\

Let now $Q\in \Psi_b(\X,\G)$ be a family of elliptic self-adjoint operators. As in \cite{a-s skew}, we consider the suspended family
\begin{equation}\label{sus fam}
Q_{sus}:= \cos(\pi t) + iQ\sin(\pi t), \quad t\in [0,1].
\end{equation}
It is an elliptic family over $B^\prime= B\times [0,1]$, invertible at $t=0, 1$, and can be regarded as an element of $\Psi_b(\X^\prime, \G)$, with $\X^\prime= \X\times [0,1]$ fibered over $B^\prime$. Moreover, $Q$ is homotopic to $Q^\prime$, through elliptic self-adjoint operators, if and only if $Q_{sus}$ is homotopic to $Q_{sus}^\prime$, through elliptic  operators.
Hence one can associate to $Q$ the class
\begin{equation}
\label{ }
[\sigma_B(Q)]_1:= [\sigma_{B^\prime}(Q_{sus})]\in K^0(T\X\times (0,1))= K^1(T\X).
\end{equation}
If $r_{\scriptsize\M}: K^1(T\X)\to K^1(T_{\scriptsize\M}\X)$ is the map of restriction to the boundary, we have that $r_{\scriptsize\M}[\sigma_B(Q)]_1$ coincides with the symbol of the indicial family of $Q_{sus}$.
In this setting, the sufficient condition for cobordism invariance given in Corollary \ref{cor susp ops} becomes
\begin{equation}
\label{odd symb susp}
r_{\scriptsize\M}[\sigma_B(Q)]_1=\beta [\sigma_B(P)].
\end{equation}
for some self-adjoint, elliptic family of $b$-pseudodifferential operators on $\X$.

The above construction applies to families of Dirac operators (see \cite{melrose APS, melrose-piazza fam dirac odd}, and also \cite{melrose-piazza fam dirac} for the odd case). Take smooth fibered manifolds $\M, \X$ over $B$, with compact, oriented  fibers $M=\partial X$, with $M$ even-dimensional. We assume that $\M$ and $\X$ are also oriented, in that the structure group reduces to orientation preserving diffeomorphisms. Assume also, as before, that $\X$ is endowed with smoothly varying $b$-metrics. Let $\E$ be an Hermitian bundle of Clifford modules over $^bT^*\X$, endowed with a smooth unitary Clifford connection on the fibres.
Let $\eth$ be the associated family of Dirac operators on $\X$; then $\eth$ is an elliptic, self-adjoint family of $b$-operators. 
Since on $\X=\M$ we have $^bT^*\X \cong \RR\left(\frac{dx}{x}\right) \oplus T^*\M$, for some boundary defining function $x$, there is an induced Clifford action of $T^*\M$ on $\E_{\scriptsize\M}$. Moreover, at the boundary, $\E$ decomposes as $\E_{\scriptsize\M}=\E^+_0\oplus \E^-_0$. 
If we denote by $\eth_0$ the induced family of Dirac operators on $\M$, then $\eth_0$ is self-adjoint, odd and elliptic.

It was noted by Melrose and Piazza in \cite{melrose-piazza fam dirac odd} that $\eth$ and $\eth_0^+$ verify (\ref{odd symb susp}) and, in fact, this condition was crucial to show that the topological family index - defined essentially through the Thom isomorphism - of  $\eth_0^+$ vanishes, so that cobordism invariance for a family of Dirac operators on a boundary then follows from the families index theorem. (See also \cite{melrose-piazza fam dirac} for the odd case and the proof). 

Now we consider families of signature operators and check that Theorem \ref{thm cbdsm inv} applies to show that the index of such a family on a boundary vanishes, using a different approach. Once we check that the relevant objects are well defined over the base, the proof in $K$-theory goes much as in the case of the single signature operator (Proposition 2.8 in \cite{cc cbdsm}).

Let $M\hookrightarrow\M\to B$ be a smooth fibration with fibre an even-dimensional, oriented manifold $M$. The bundle $\Lambda^*(T^*\M)$ of forms along the fibres is a smooth bundle over $\M$ with fibre diffeomorphic to $\Lambda(T^*M)$. 
 There is a smooth action of $\Lambda^*(T^*\M)$ on itself,
 \begin{equation}
\label{action cliff}
\Lambda^*(T^*\M) \times \Lambda^*(T^*\M) \to \Lambda^*(T^*\M),\quad (\xi, \eta)\mapsto c(\xi)\eta,
\end{equation}
such that for $\xi \in T^*\M$, $c(\xi)\eta= \xi \wedge \eta - i_\xi(\eta)$. 
We also have a $\mathbb{Z}_2$-grading $\Lambda^*(T^*\M)=\Lambda^+(T^*\M)\oplus \Lambda^-(T^*\M)$, where $\Lambda^\pm(T^*\M):= (1\pm \Gamma)\Lambda^*(T^*\M)$, with $\Gamma:= i^{n\slash 2}c(e_1)...c(e_n)$ and $e_1,...,e_n$ an oriented basis for $T_x\M$, $x\in \M$. In particular, $c(\xi) \Lambda^\pm(T^*\M)\subset \Lambda^\mp(T^*\M)$, for all $\xi\in T^*\M$. 

 Let $d$ denote the exterior derivative along the fibres and $d^*$ its adjoint. We define the \textit{signature family} on $\M$ as
\begin{equation}
\label{sign fam}
D=d+d^*: C^\infty(\M; \Lambda^+(T^*\M)) \to C^\infty(\M; \Lambda^-(T^*\M)).
\end{equation}
Clearly, $D$ restricts to signature operators on the fibres. The symbol of $D$ is given by the map $\sigma(D): \pi^* \Lambda^+(T^*\M) \to \pi^* \Lambda^-(T^*\M)$, with $\pi: T^*\M \to \M$ the projection, such that $\sigma(D)(\xi)= c(\xi)$. Hence, $\sigma(D)$ is invertible outside the zero-section, that is, $D$ is elliptic. If $\M$ is compact, then it defines a $K$-theory class
\begin{equation}
\label{symb sign fam}
 [\sigma(D)]=[\pi^*\Lambda^+(T^*\M), \pi^*\Lambda^-(T^*\M), c]\in K^0(T^*\M).
\end{equation}
 Note that in this case $\ker D_b$ and $\coker D_b$ are vector bundles over $B$, hence the \textit{signature index class} is given by
\begin{equation}
\label{ind sign fam}
\ind (D)= [\ker D] -[\coker D]\in K^0(B).
\end{equation}
We check that it is an invariant of cobordism:

\begin{cor}
Let $D$ be the signature family on a fibered manifold $M\hookrightarrow \M\to B$, with $M$ an even-dimensional, oriented compact manifold. If there is a compact, oriented fibered manifold $X\hookrightarrow\X\to B$ with $\partial \X=\M$, then $\ind (D)=0$.
\end{cor}

\pf 
The proof follows the lines of the proof of Proposition 2.8 in \cite{cc cbdsm}. We sketch the main points. First consider the signature family $D_1$ on $X\times\RR\hookrightarrow \X\times \RR\to B$. Restricting its symbol to $\X\times\{0\}$, we get a $K$-theory class
$$\omega:=  [(\sigma(D_1)_{|\scriptsize{T\X\times \RR}}]\in K^0(T\X\times\RR)=K^1(T\X).$$
We have that $r(\omega)= [(\sigma(D_2))_{|\scriptsize{\M\times\{0\}}}]$, where $r: K^1(T\X)\to K^0(T\M \times \RR^2) $ is the restriction map and $D_2$ is the signature family on $\M\times\RR^2$. 
Now, if $\beta$ denotes a representative for the Bott class in $K^0(\RR^2)$, then $\sigma(D_2)= \sigma(D)\cup(\pi^*{\beta}\oplus\pi^*{\beta})$, so that 
$$r(\omega)=(\sigma(D)\cup\beta)\oplus(\sigma(D)\cup\beta)=\beta_{T\scriptsize{\M}}(\sigma(D)\oplus\sigma(D)).$$
We conclude that $(\M, \sigma(D)\oplus\sigma(D))\sim 0$. It follows from Theorem \ref{thm cbdsm inv} that $\ind(D\oplus D)=0$. From the additivity of the family index, we have then $\ind(D)=\ind(D^*)\Leftrightarrow [\ker D] -[\coker D]=[\ker D^*] -[\coker D^*]$, so that there exist trivial bundles $\theta^m$ and $\theta^n$ such that $\ker D\oplus \theta^n\cong \ker D^*\oplus \theta^m$ and $\coker D\oplus \theta^n\cong \coker D^*\oplus \theta^m$. Since $\ker D^*\cong \coker D$, we have $n=m$ and therefore $ [\ker D] =[\coker D]$. Hence, $\ind(D)=0$.

\qed

A similar result holds for signature families twisted by some smooth bundle $\W$ over $\M$ with a smooth fibre connection, as long as $\W$ can be extended to the boundary.

To finish, we give an alternative $K$-theory formulation of cobordism invariance, following Moroianu in \cite{moroianu cc cbdsm}.
 By a metric on $\M$ we mean a continuous family of smooth metrics on $TM_b$, $b\in B$.  The unit sphere bundle $S^*\M=S(T^*\M)$ is well defined and it is a manifold over $B$. We let $T^*_{sus}\M:=T^*\M\times \RR$, as before, and $S^*_{sus}\M=S(T^*_{sus}\M)$. Then, as Moroianu noted, there is an isomorphism
\begin{equation}
\label{ }
d: K^0(T^*\M)\to K^0(S^*_{sus}\M)\slash \pi^*K^0(\M),
\end{equation}
where $\pi: S^*_{sus}\M\to \M$ is the projection, such that given vector bundles $\E^+$, $\E^-$ over $T^*\M$ (not necessarily smooth on fibers) and a map $\sigma: \E^+ \to \E^-$, which is an isomorphism outside the unit ball, we have
\begin{equation}
\label{ }
\quad [\E^+,\E^-,\sigma]\in K^0(T^*\M)\mapsto 
\begin{cases}
\E^+,& \text{ on } S^*\M\cap \{\xi\geq 0\},\\
\E^-,& \text{ on } S^*\M\cap \{\xi< 0\},
\end{cases}
\end{equation}
with $\E^+$, $\E^-$ identified via $\sigma$ on $S^*_{sus}\M\cap \{\xi= 0\}=S^*\M$.
Moreover, taking the boundary maps for the relative pairs $(B^*\X, S^*\X)$ and $(B^*_{sus}\M, S^*_{sus}\M)$, where $\M=\partial_B \X$, and the maps of restriction to the boundary, the following diagram commutes:
\[
\begin{CD}
K^0(S^*\X)                           @>{\partial}>>        K^1(T^*\X)     \\
@V{r}VV                                     @VV{r}V           \\
K^0(S^*_{sus}\M)         @>{\partial}>>           K^1(T^*_{sus}\M)\\  
@V{q}VV                                     @VV{\beta^{-1}}V           \\
K^0(S^*_{sus}\M) \slash \pi^*K^0(\M)        @>>{d^{-1}}>           K^0(T^*\M).
\end{CD}
\]
Since $\beta^{-1}\circ r =u$, the map of restriction of symbols (\ref{map rest symb}), we have the following analogue of Moroianu's $K$-theory formulation of cobordism invariance for families.

\begin{thm}\label {thm cbdsm inv moroianu}
Let $\M$, $\X$ be manifolds over $B$ with $\M=\partial_B\X$, and $P$ be an elliptic family of pseudodifferential operators over $B$ with principal symbol $\sigma\in K^0(T^*\M)$. If $d(\sigma)\in r(K^0(S^*\X) ) +  \pi^*K^0(\M) $, then $\ind(P)=0$. 
\end{thm}

\pf If $d(\sigma)= r(\omega) +  \pi^*K^0(\M) $, with  $\omega \in K^0(S^*\X) $, then $u(\partial \omega)=\sigma$; hence $(\M,\sigma)\sim 0$. 
\qed\\

\noindent\textit{Acknowledgments}. 

 I would like to thank Victor Nistor and Sergiu Moroianu for very helpful discussions. I would also like to thank  Michel Hilsum and Radu Popescu and the referees for useful comments which helped improve the paper.

\end{document}